\documentclass[12pt]{amsart}
\usepackage{verbatim}
\usepackage{amssymb,amsmath}

\newtheorem{lemma}{Lemma} [section]
\newtheorem{thm}[lemma]{Theorem}
\newtheorem{cor}[lemma]{Corollary}

\theoremstyle{remark}
\newtheorem*{remark}{Remark}

\newcommand{\rmv}[1]{}
\newcommand{\F}{{\mathbb F}}

\DeclareMathOperator{\chara}{char}
\newcommand{\Z}{{\mathbb Z}}

\newcommand{\C}{{\mathbb C}}

\begin{document}



\title{Rational functions with linear relations}

\author{Ariane M. Masuda}
\address{
  School of Mathematics and Statistics,
  Carleton University,
  1125 Colonel By Drive,
  Ottawa, ON K1S 5B6,
  Canada
}
\email{ariane@math.carleton.ca}

\author{\hbox{Michael E. Zieve}}
\address{
  Center for Communications Research
  805 Bunn Drive
  Princeton, NJ 08540
}
\email{zieve@math.rutgers.edu}

\thanks{The authors thank Bob Beals, Alan Beardon, Alex Eremenko
and Patrick Ng for useful correspondence.}



\begin{abstract}
We find all polynomials $f,g,h$ over a field $K$ such that $g$ and
$h$ are linear and $f(g(x))=h(f(x))$.  We also solve the same
problem for rational functions $f,g,h$, in case the field $K$ is
algebraically closed.
\end{abstract}

\date{\today}

\maketitle


\section{Introduction}

Around 1920, Fatou, Julia and Ritt made profound investigations of
functional equations.  In particular, they wrote at length on
commuting rational functions: that is, $f,g\in\C(x)$ with $f(g(x))=g(f(x))$.
Fatou and Julia \cite{Fatoucomm, Juliacomm} found all solutions when the
Julia set of $f$ or $g$ is not the Riemann sphere.  This includes the case
of polynomials of degree at least $2$, where up to conjugacy by a linear
polynomial, either $f=x^n$ and $g=x^m$ are
power polynomials, or $f=T_n$ and $g=T_m$ are Chebychev polynomials,
or $f$ and $g$ have a common iterate.  Using different methods which did
not require the Julia set hypothesis, Ritt \cite{Ritt20} determined
precisely when two polynomials have a common
iterate, and moreover \cite{Ritt23} he found all commuting rational
functions.  Years later, Eremenko \cite{Eremenko} proved Ritt's results
using methods of modern iteration theory.

Julia showed that commuting rational functions have the same Julia set.
Conversely, much subsequent work has shown that rational functions with
the same Julia set are related to commuting rational functions (cf.\
\cite{LP} and the references therein).  In particular, for polynomials
this relationship involves composition with a rotational symmetry of
the Julia set.

Several authors have considered analogous questions over fields $K$ of
positive characteristic, but there are few satisfactory results.
There are new types of examples, for instance any two additive
polynomials $\sum_i a_i x^{p^i}$ over the prime field $\F_p$ commute.

In fact, challenges arise already in finding the commuting polynomials
$f,g\in K[x]$ in the special case $\deg(g)=1$.
Wells \cite{Wells} and Mullen \cite{Mullen} solved this problem
over finite fields $K$, so long as $\deg(f)<\#K$.  Park \cite{Park} proved
similar results.
Eigenthaler and N\"obauer \cite{EN} solved the problem in various
special cases, for instance if $\deg(f)=\chara(K)$.  In this paper we
solve the problem in general, and more generally we find all $f,g,h\in K[x]$
with $\deg(g)=\deg(h)=1$ such that $f\circ g=h\circ f$:

\begin{thm}
\label{thmpol}
Let $K$ be a field.  The entries in the following list with $f\notin K$
comprise all values $\alpha,\beta,\gamma,\delta\in K$ and $f\in K[x]\setminus K$ such that
$\alpha,\gamma\ne 0$ and $f(\alpha x+\beta) = \gamma f(x)+\delta$:
\begin{enumerate}
\item[(0)] $\alpha =\gamma =1$, $\beta=\delta=0$, and $f\in K[x]$;
\item[(1)] $\alpha =\gamma =1$, $\beta\ne 0$, and $f=(\delta/\beta)x+r$ with
\[
\begin{cases}
r\in K & \text{if $\chara(K)=0$} \\
r\in K[x^p-\beta^{p-1}x] & \text{if $\chara(K)=p>0$};
\end{cases}\]
\item[(2)] $\alpha\ne 1$, $\gamma=\alpha^e$, and $f=f_0+f_1(x-\beta/(1-\alpha))$, where
$e,s\in\Z_{\ge 0}$, $f_1\in x^e K[x^s]$, $\alpha^s=1$,
and $f_0\in K$ satisfies $\delta=(1-\gamma)\cdot f_0$.
\end{enumerate}
\end{thm}

In case $K=\C$, the polynomials in (2) are those for which
the Julia set has a rotational symmetry \cite{BE}.  Moreover, Ritt
showed \cite{Ritt22}
that the decomposition of a complex polynomial into indecomposables is unique,
except for nonuniqueness coming from composing a linear with its inverse,
or using the commutativity of Chebychev polynomials, or using the identity
$x^s \circ x^e \psi(x^s) = x^e \psi(x)^s \circ x^s$.  Note that the polynomials
$x^e\psi(x^s)$ from this identity occur in (2).  Ritt's identity has
a characteristic~$p$ analogue \cite{BZ}, namely
 $(x^p-x)\circ (x+\psi(x^p-x)) =
(x + \psi^p-\psi)\circ (x^p-x)$, and it is interesting that the
polynomials $x+\psi(x^p-x)$ occur in (1).

We also prove an analogous result for rational functions:

\begin{thm}
\label{thmrat}
Let $K$ be a field of characteristic $p\ge 0$.
The entries in the following list with $f\notin K$ comprise all
$g,h\in K(x)$ and $f\in K(x)\setminus K$ such that $f\circ g=h\circ f$
and each of $g$ and $h$ has degree one and
has a fixed point in $K\cup\{\infty\}$; here $u,v,\psi\in K(x)$
and $\deg(u)=\deg(v)=1$\emph{:}
\begin{enumerate}
\item $f=u^{-1}\circ (\delta x+\psi(x^p-x))\circ v^{-1}$, $g=v(v^{-1}(x)+1)$,
and $h=u^{-1}(u(x)+\delta)$, where $\delta\in K$, and
if $p=0$ then $\psi\in K$;
\item $f=u^{-1}\circ x^e \psi(x^s)\circ v^{-1}$, $g=v(\alpha v^{-1}(x))$,
$h=u^{-1}(\alpha^e u(x))$, where $e,s\in\Z$, $\alpha\in K^*$, and
$\alpha^s=1$.
\end{enumerate}
\end{thm}

Our hypothesis on fixed points is always true if $K$ is algebraically
closed.  To apply this result to arbitrary fields $K$, one might need
$u,v,\psi$ to have coefficients in an extension of $K$.

In case $K=\C$, these results were proved by af H\"allstr\"om
\cite{Halbpol,Halbrat}.  His method has some features in common with ours,
but is somewhat more complicated.

We give a quick inductive proof of Theorem~\ref{thmpol} in the next
section.  Then in Sections~\ref{sec scale} and \ref{sec rat} we use ideas
from dynamics and Galois theory to prove Theorem~\ref{thmrat}, which
yields another proof of Theorem~\ref{thmpol}.  Finally, in
Section~\ref{sec mullen} we deduce the results of Wells~\cite{Wells},
Mullen~\cite{Mullen} and Park~\cite{Park}
as consequences of Theorem~\ref{thmpol}.


\section{Polynomial solutions}

In this section we prove Theorem~\ref{thmpol}.

Pick $\alpha,\beta,\gamma,\delta\in K$ with $\alpha,\gamma\ne 0$, and let 
$f\in K[x]$ have degree $n>0$.
We will determine when $f(\alpha x+\beta)=\gamma f(x)+\delta$.
We assume $\gamma=\alpha^n$, since otherwise $f(\alpha x+\beta)$ and 
$\gamma f(x)+\delta$ have distinct leading coefficients.
First suppose $\alpha=1$, so $\gamma=1$.  If $\beta=0$ then our equation
becomes $f(x)=f(x)+\delta$, 
so $\delta=0$; conversely, if $\beta=\delta=0$ then trivially every $f$ is a 
solution.  So assume $\beta\ne 0$, and put $r:=f-(\delta/\beta)x$; then $f$ 
satisfies $f(x+\beta)=f(x)+\delta$ if and only if $r$ satisfies $r(x+\beta)=r(x)$.  
Let $p:=\mathrm{char}(K)$ and $m:=\deg(r)$.  If
$p\nmid m$ then there are no such $r$, since $r(x+\beta)-r(x)$ has degree $m-1$.
In particular, if $p=0$ then $r\in K$, so assume $p>0$.
Plainly every $\hat r\in K[x^p-\beta^{p-1}x]$ satisfies $\hat r(x+\beta)=
\hat r(x)$.  For any $r\in K[x]$ with $r(x+\beta)=r(x)$, we know that
$p\mid \deg(r)$, so there is some $\hat r\in K[x^p-\beta^{p-1}x]$
which has the same leading term as $r$; but then $\tilde r:=r-\hat r$ satisfies
$\tilde r(x)=\tilde r(x+\beta)$ and $\deg(\tilde r)<\deg(r)$, so it follows by
induction on $\deg(r)$ that $r\in K[x^p-\beta^{p-1}x]$.

Now suppose $\alpha\ne 1$.  Let $s$ be the multiplicative order of $\alpha$, if this
order is finite; otherwise put $s=0$.  Thus the integers $m$ with $\alpha^m=1$
are precisely the multiples of $s$.  Let $u$ be the leading coefficient of
$f$, and put $\hat f:=u\cdot (x - \beta/(1-\alpha))^n$; then
\[
\hat f(\alpha x+\beta) = u\cdot (\alpha x-\alpha\beta/(1-\alpha))^n = 
u\alpha^n \cdot (x-\beta/(1-\alpha))^n = \gamma\hat f(x).
\]
Now put $\tilde f:=f-\hat f$, and note that $\tilde n:=\deg(\tilde f)<n$;
moreover, $f$ satisfies $f(\alpha x+\beta)=\gamma f(x)+\delta$ if and only if 
$\tilde f$ satisfies $\tilde f(\alpha x+\beta) = \gamma\tilde f(x) +\delta$.
If $\tilde n>0$, then the leading coefficients of $\tilde f(\alpha x+\beta)$ and
$\gamma\tilde f(x)+\delta$ are identical if and only if $\alpha^{\tilde n}=\gamma=\alpha^n$, or
equivalently $\tilde n\equiv n\pmod{s}$.
By induction on $\deg(f)$, it follows that $f$ satisfies $f(\alpha x+\beta)=\gamma f(x)+\delta$
if and only if $f=f_0+f_1(x-\beta/(1-\alpha))$ where $f_0\in K$ satisfies
$f_0=\gamma f_0+\delta$ and $f_1\in xK[x]$ has only
terms of degree congruent to $n\pmod{s}$.  The result follows.
\hfill\qed


\section{Solutions involving scalings or translations}
\label{sec scale}

In this section we solve the equation $f\circ g=h\circ f$ in rational
functions $f,g,h\in K(x)$ where $g,h\in xK^*\cup\{x+1\}$.
Here $K$ is a field of characteristic $p\ge 0$.  Let $L=K(x^p-x)$ if $p>0$,
and put $L=K$ if $p=0$.

\begin{lemma}
\label{translate0}
For $f\in K(x)$, we have $f(x+1)=f(x)$ if and only if $f\in L$.
\end{lemma}

\begin{proof}
Let $\sigma$ be the $K$-automorphism of $K(x)$ which maps $x\mapsto x+1$.
Then $L$ is the subfield of $K(x)$ fixed by $\sigma$.
Thus $f\in L$ if and only if $\sigma(f)=f$, or equivalently $f(x+1)=f(x)$.
\end{proof}

\begin{cor}
\label{translatelr}
For $f\in K(x)$, we have $f(x+1)=f(x)+1$ if and only if $f-x\in L$.
\end{cor}

\begin{proof}
Putting $r(x):=f(x)-x$, we have $f(x+1)=f(x)+1$ if and only if
$r(x+1)=r(x)$, so the result follows from Lemma~\ref{translate0}.
\end{proof}

\begin{cor}
\label{translater} For any $\gamma\in K$ and $f\in K(x)$ with $f\ne 0$,
we have $f(x+1)=\gamma f(x)$ if and only if $\gamma=1$ and $f\in L$.
\end{cor}

\begin{proof}
The leading terms of both the numerator and denominator of $f(x)$
are identical to those of $f(x+1)$, so if $f(x+1)=\gamma f(x)$ then $\gamma=1$;
now the result follows from Lemma~\ref{translate0}.
\end{proof}

\begin{lemma}
\label{scalepoly} For any $\alpha,\gamma\in K^*$ and any nonzero $f\in K[x]$,
we have $f(\alpha x)=\gamma f(x)$ if and only if $f=x^e \psi(x^s)$ for some
$\psi\in K[x]$ and $e,s\in\Z_{\ge 0}$ with $\alpha^e=\gamma$ and $\alpha^s=1$.
\end{lemma}

\begin{proof}
Equate coefficients of terms of the same degrees in $f(\alpha x)$ and
$\gamma f(x)$.
\end{proof}

\begin{cor}
\label{scale} For any $\alpha,\gamma\in K^*$ and any nonzero $f\in K(x)$, we
have $f(\alpha x)=\gamma f(x)$ if and only if $f=x^e \psi(x^s)$ for some
$\psi\in K(x)$ and $e,s\in\Z$ with $\alpha^e=\gamma$ and $\alpha^s=1$.
\end{cor}

\begin{proof}
The `if' direction is clear, so suppose $f(\alpha x)=\gamma f(x)$. Write
$f=f_1/f_2$ with coprime $f_1,f_2\in K[x]$.  Then the denominators
of $f(\alpha x)$ and $\gamma f(x)$ are $f_2(\alpha x)$ and $f_2(x)$, so $f_2(\alpha
x)=\eta f_2(x)$ for some $\eta\in K^*$.  Thus $f_1(\alpha x)=\gamma\eta
f_1(x)$.  Now the result follows by applying Lemma~\ref{scalepoly}
to both $f_1$ and $f_2$.
\end{proof}

\begin{lemma}
\label{translatel}
For $\alpha\in K^*$ and $f\in K(x)$, we have $f(\alpha x)\ne f(x)+1$.
\end{lemma}

\begin{proof}
Suppose to the contrary that $f(\alpha x)=f(x)+1$.  Plainly $x=0$ must be
a pole of $f$.  Write $f=f_1/f_2$ with coprime $f_1,f_2\in K[x]$ (so
$f_1(0)\ne 0$ and $f_2(0)=0$, whence $\deg(f_2)>0$).  Then the
denominators of $f(\alpha x)$ and $f(x)+1$ are $f_2(\alpha x)$ and $f_2(x)$,
so we have $f_2(\alpha x)=\eta f_2(x)$ for some $\eta\in K^*$.  Then
\[
\frac{f_1(\alpha x)}{f_2(\alpha x)} = \frac{f_1(x)}{f_2(x)} + 1
\]
implies that
\[
\frac{f_1(\alpha x)}{\eta} = f_1(x) + f_2(x).
\]
Since $f_1(0)\ne 0$ and $f_2(0)=0$, substituting $x=0$ gives $\eta=1$.
Thus $f_2(\alpha x)=f_2(x)$; since $f_2$ is nonconstant, it follows that
$s:=\#\langle\alpha\rangle <\infty$ and $f_2\in K[x^s]$.  But then
$f_1(\alpha x) - f_1(x) = f_2(x)\in K[x^s]$, which is impossible since
$f_1(\alpha x)-f_1(x)$ has no terms of degree divisible by~$s$.
\end{proof}


\section{Solutions with arbitrary linears}
\label{sec rat}

In this section we solve the equation $f\circ g = h\circ f$ in rational
functions $f,g,h\in K(x)$ with $\deg(g)=\deg(h)=1$.  Here $K$ is a field
of characteristic $p\ge 0$.  We will reduce to the
cases considered in the previous section, by means of the following
observation: if $u,v\in K(x)$ satisfy $\deg(u)=\deg(v)=1$, then $f\circ g =
h\circ f$ if and only if $F\circ G = H\circ F$, where $F:=u\circ f\circ v$,
$G:=v^{-1}\circ g\circ v$, and $H:=u\circ h\circ u^{-1}$.

First consider the case of polynomials $f,g,h\in K[x]$.  Then $g=\alpha x+\beta$ for
some $\alpha,\beta\in K$ with $\alpha\ne 0$.  If $\alpha\ne 1$ then 
$v:=x+\beta/(1-\alpha)$ satisfies $v^{-1}\circ g\circ v=\alpha x$.  If $\alpha=1$ and 
$\beta\ne 0$ then $v:=\beta x$ satisfies
$v^{-1}\circ g\circ v=x+1$.  Thus, in any case there is a degree-one
$v\in K[x]$ such that $G:=v^{-1}\circ g\circ v$ is either $\alpha x$ or $x+1$.
Likewise, writing $h=\gamma x+\delta$, there is a degree-one $u\in K[x]$ such that
$H:=u\circ h\circ u^{-1}$ is either $\gamma x$ or $x+1$.  Now the above observation,
in combination with the results of the previous section, implies the following
version of Theorem~\ref{thmpol}:

\begin{thm}
The polynomials $f,g,h\in K[x]$ such that $f\circ g=h\circ f$ and
$\deg(g)=\deg(h)=1\le\deg(f)$ are as follows; here
$u,v,\psi\in K[x]$ and $\deg(u)=\deg(v)=1$\emph{:}
\begin{enumerate}
\item $f=u^{-1}\circ (x+\psi(x^p-x))\circ v^{-1}$, $g=v(v^{-1}(x)+1)$,
and $h=u^{-1}(u(x)+1)$, where if $p=0$ then $\psi\in K$;
\item $p>0$, $f=\psi(x^p-x)\circ v^{-1}$, $g=v(v^{-1}(x)+1)$, and $h=x$,
where $\deg(\psi)>0$;
\item $f=u^{-1}\circ x^e \psi(x^s)\circ v^{-1}$, $g=v(\alpha v^{-1}(x))$,
$h=u^{-1}(\alpha^e u(x))$, where $e,s\in\Z_{\ge 0}$, $\alpha\in K^*$,
$\alpha^s=1$, and $\deg(x^e\psi(x^s))>0$.
\end{enumerate}
\end{thm}

Note that we can combine the first two possibilities into the single
possibility $f=u^{-1}\circ (\delta x+\psi(x^p-x))\circ v^{-1}$,
$g=v(v^{-1}(x)+1)$, and $h=u^{-1}(u(x)+\delta)$ with $\delta\in K$.

Next we consider rational functions.  Any degree-one $g\in K(x)$ has a fixed
point $\rho$, though this fixed point might lie in a quadratic extension of
$K$.  If $\rho=\infty$ then $g\in K[x]$; if $\rho\ne\infty$ then for
$v=\rho+1/x$ we see that $v^{-1}\circ g\circ v$ fixes $\infty$, and hence
lies in $K[x]$.  We can then proceed as above, resulting in a proof of
Theorem~\ref{thmrat}.


\section{Derivation of prior results}
\label{sec mullen}

In this section we explain how Theorem~\ref{thmpol} relates to the
results of Wells~\cite{Wells}, Mullen~\cite{Mullen}, and
Park~\cite{Park}, all of which were formulated in a quite different
manner.  

It follows from Theorem~\ref{thmpol} that, if $K$ is a
field of characteristic $p>0$, and if we prescribe elements $\beta,\delta\in K$
with $\beta\ne 0$, then the polynomials $f\in K[x]$ such that
$f(x+\beta)=f(x)+\delta$ are precisely the elements of $(\delta/\beta)x+K[x^p-\beta^{p-1}x]$.
In particular, writing $f=\sum_{i=0}^n f_i x^i$, we see that the
coefficients $f_0,f_p,f_{2p},\dots$ can be arbitrary elements of $K$,
and these coefficients uniquely determine all the other $f_i$'s.
This generalizes the results of Wells and Park.

Wells~\cite{Wells} restricted to the case that $K=\F_q$ is finite, $\deg(f)<q$,
and $\delta=\beta$.  Wells used a different method.  Namely, by considering
terms of degree $x^{i-1}$ in the functional equation $f(x+\beta)=f(x)+\beta$,
one can solve for $if_i$ in terms of the coefficients $f_j$ with
$j>i$; thus, by successively computing $f_n, f_{n-1},\dots,f_1$, we
see that the coefficients $f_i$ with $p\nmid i$ are uniquely
determined by the coefficients $f_{pj}$.  Hence there are at most
$q^{q/p}$ possibilities for $f$; but this equals the number of
mappings $\F_q\to \F_q$ which commute with the map $x\mapsto x+\beta$.
Since every mapping $\F_q\to \F_q$ is induced by a unique polynomial
of degree less than $q$, it follows that the $f_{pj}$ can be
arbitrary elements of $\F_q$.  On the other hand, as noted above,
this fact follows at once from our expression
$x+K[x^p-\beta^{p-1}x]$ for all such $f$'s.

Park~\cite{Park} considered the case that $K$ is finite,
$\deg(f)<p^2$, and $\beta,\delta\in K^*$.  He wrote out the conditions
on the $f_j$'s coming from equating terms of like degrees in the
functional equation $f(x+\beta)=f(x)+\delta$, and proved his result via
several pages of calculations involving binomial coefficients.
In these calculations, the hypothesis $\deg(f)<p^2$ yielded
crucial simplifications.

Next suppose $\alpha\in K^*$ is a primitive $s^{\operatorname{th}}$
root of unity, and suppose $\gamma=\alpha^e$ with $0<e<s$.  Fix $\beta,\delta\in K$.
By Theorem~\ref{thmpol}, the polynomials $f\in K[x]$ such that
$f(\alpha x+\beta)=\gamma f(x)+\delta$ are precisely the elements of
\[ \frac{\delta}{1-\gamma} + \left(x-\frac{\beta}{1-\alpha}\right)^e
  K\left[\left(x-\frac{\beta}{1-\alpha}\right)^s\right].\]
In particular, writing $f=\sum_{i=0}^n f_i x^i$, the coefficients
$f_{e+sj}$ can be arbitrary elements of $K$, and these
coefficients uniquely determine all the other $f_i$'s.
This generalizes the main result proved by Mullen.

Mullen~\cite{Mullen} restricted to the case that $K=\F_q$ is finite,
$\deg(f)<q$,
$\gamma=\alpha$ and $\delta=\beta$.  He used the same method as Wells: equating
coefficients of $x^i$ in $f(\alpha x+\beta)=\alpha f(x)+\beta$ enables one to express
$(\alpha^i-\alpha)f_i$ in terms of $f_{i+1},f_{i+2},\dots,f_n$.
Since $\alpha^i\ne\alpha$ if $i\not\equiv 1\pmod{s}$, it follows that all
the $f_i$'s are uniquely determined by the coefficients $f_{1+sj}$.
Hence there are at most $q^{(q-1)/s}$ possibilities for $f$; but
this equals the number of mappings $\F_q\to\F_q$ which commute
with $x\mapsto \alpha x+\beta$, so the coefficients $f_{1+sj}$ can be
arbitrary elements of $\F_q$.

\begin{remark}
As stated, \cite[Thm.~1]{Mullen} asserts
that two polynomials are equal if they have the same coefficients.
The proof in \cite{Mullen} shows that the result would remain true
(and become nontrivial) if we require $p\nmid s$ when $b=1$, and
$\#\langle b\rangle\nmid (s-1)$ if $b\ne 1$.  Our comments above
refer to this corrected version of Mullen's result.  Also, the
papers \cite{Wells,Mullen} comment on polynomials over $\F_q$ of
degree $\ge q$ which commute with linear polynomials, but in
those papers commutation is only studied modulo $x^q-x$; in other words,
they consider $f(\alpha x+\beta)=\alpha f(x)+\beta$ as an equality of
functions on $\F_q$, rather than an equality of polynomials.
\end{remark}

To summarize, it seems that the complications Wells, Mullen and Park
encountered were caused by their desire to phrase the results in terms
of the coefficients of $f$ as an element of $K[x]$; the key to our
simpler presentation is that we directly represent $f$ in terms of an
additive subgroup of $K[x]$.



\begin{thebibliography}{9}
\newcommand{\au}[1]{{#1},}
\newcommand{\ti}[1]{\textit{#1},}
\newcommand{\jo}[1]{{#1}}
\newcommand{\vo}[1]{\textbf{#1}}
\newcommand{\yr}[1]{(#1),}
\newcommand{\pp}[1]{#1.}
\newcommand{\pps}[1]{#1;}
\newcommand{\bk}[1]{{#1},}
\newcommand{\inbk}[1]{in: {#1}}
\newcommand{\xxx}[1]{{arXiv:#1}}

\bibitem{Halbpol}
\au{G. af H\"allstr\"om}
\ti{\"Uber halbvertauschbare Polynome}
\jo{Acta Acad. Abo.}
\vo{21}
\yr{1957}
no.\ 2, 20 pp.

\bibitem{Halbrat}
\au{G. af H\"allstr\"om}
\ti{\"Uber Halbvertauschbarkeit zwischen linearen und allgemeineren
rationalen Funktionen}
\jo{Math. Japon.}
\vo{4}
\yr{1957}
\pp{107--112}

\bibitem{BE}
\au{I.~N. Baker and A. Er\"emenko}
\ti{A problem on Julia sets}
\jo{Ann. Acad. Sci. Fenn.}
\vo{12}
\yr{1987}
\pp{229--236}

\bibitem{BZ}
\au{R.~M. Beals and M.~E. Zieve}
\ti{Decompositions of polynomials}
preprint, 2007.

\bibitem{EN}
\au{G. Eigenthaler and W. N\"obauer}
\ti{\"Uber die mit einem Polynom vertauschbaren linearen Polynome}
\jo{\"Osterreich. Akad. Wiss. Math.-Natur. Kl. Sitzungsber. II}
\vo{199}
\yr{1990}
\pp{143--153}

\bibitem{Eremenko}
\au{A.~\`E. Er\"emenko}
\ti{Some functional equations connected with
the iteration of rational functions}
\jo{Algebra i Analiz}
\vo{1}
\yr{1989}
\pp{102--116}
(Translated in \jo{Leningrad Math. J.}
\vo{1}
\yr{1990}
\pp{905--919})

\bibitem{Fatoucomm}
\au{P. Fatou}
\ti{Sur l'iteration analytique et les substitutions permutables}
\jo{J. Math. Pures Appl. (9)}
\vo{2}
\yr{1923}
\pp{343--384}

\bibitem{Juliacomm}
\au{G. Julia}
\ti{M\'emoire sur la permutabilit\'e des fractions rationnelles}
\jo{Ann. Acad. \'Ecole Norm. Sup.}
\vo{39}
\yr{1922}
\pp{131--215}

\bibitem{LP}
\au{G.~M. Levin and F. Przytycki}
\ti{When do two rational functions have the same Julia set?}
\jo{Proc. Amer. Math. Soc.}
\vo{125}
\yr{1997}
\pp{2179--2190}

\bibitem{Mullen}
\au{G.~L. Mullen} \ti{Polynomials over finite fields which commute
with linear permutations} \jo{Proc. Amer. Math. Soc.} \vo{84}
\yr{1982} \pp{315--317}

\bibitem{Park}
\au{H.~G. Park} \ti{Polynomials satisfying $f(x+a)=f(x)+c$ over
finite fields} \jo{Bull. Korean Math. Soc.} \vo{29} \yr{1992}
\pp{277--283}

\bibitem{Ritt20}
\au{J.~F. Ritt}
\ti{On the iteration of rational functions}
\jo{Trans. Amer. Math. Soc.}
\vo{21}
\yr{1920}
\pp{348--356}

\bibitem{Ritt22}
\au{J.~F. Ritt}
\ti{Prime and composite polynomials}
\jo{Trans. Amer. Math. Soc.}
\vo{23}
\yr{1922}
\pp{51--66}

\bibitem{Ritt23}
\au{J.~F. Ritt}
\ti{Permutable rational functions}
\jo{Trans. Amer. Math. Soc.}
\vo{25}
\yr{1923}
\pp{399--448}

\bibitem{Wells}
\au{C. Wells} \ti{Polynomials over finite fields which commute with
translations} \jo{Proc. Amer. Math. Soc.} \vo{46} \yr{1974}
\pp{347--350}

\end{thebibliography}
\end{document}